\magnification=\magstep1
\input amstex
\documentstyle{amsppt}
\hsize = 6.5 truein
\vsize = 9 truein

\UseAMSsymbols

\loadbold

\define \J {\Bbb J}
\define \rL {\Bbb L}
\define \K {\Bbb K}
\define \rO {\Bbb O}
\define \N {\Bbb N}
\define \AN {\Bbb N^*}
\define \rP {\Bbb P}
\define \CP {\Cal P}
\define \CET {\Cal C (ET)}
\define \s {\Bbb S}
\define \zs {\Bbb S_0}
\define \CRET {\Cal C (RET)}
\define \Z {\Bbb Z}
\define \Zi {\Bbb Z[i]}
\define \Q {\Bbb Q}
\define \R {\Bbb R}
\define \C {\Bbb C}
\define \oN {\overline {\Bbb N}}
\define \CPET {\Cal C (PET)}
\define \CePET {\Cal C_e (PET)}
\define \oR {\overline {\Bbb R}}
\define \RP {\R^{+*}}
\define \cal {\text {cal}}
\define \th	{\theta}
\define \cp {\text {cap}}
\define \gri {\text {gri}}
\define \lgr {\text {lgr}}
\define \grs {\text {grs}}
\define \ugr {\text {ugr}}
\define \gr {\text {gr}}
\define \e	{\epsilon}
\define \de	{\delta}
\define \p {\Bbb P}
\define \Gi {\Bbb Z[i]}

\define \ee {\'e}
\define \ea {\`e}
\define \oo {\^o}
\define \ž  {\^u}
\define \" {\^\i}
\define \‰ {\^a}

\define \CC {\Cal C}
\define \CF {\Cal F}
\define \CG {\Cal G}
\define \CH {\Cal H}
\define \CI {\Cal I}
\define \CN {\Cal N}
\define \CT {\Cal T}
\define \CU {\Cal U}
\define \CV {\Cal V}
\define \CW {\Cal W}
\define \CA {\Cal A}
\define \CB {\Cal B}
\define \CE {\Cal E}
\define \CR {\Cal R}
\define \CS {\Cal S}
\define \CM {\Cal M}
\define \CO {\Cal O}
\define \CX {\Cal X}

\define \card {\text{card}}

\define \bS {\bf S \rm}

\define \noo {$\text{n}^o$}

\define \T- {\overset {-1}\to T}
\define \G- {\overset {-1}\to G}
\define \GA- {\overset {-1}\to \Gamma}
\define \RA- {\overset {-1}\to R}
\define \PA- {\overset {-1}\to P}
\define \ig {g^{-1}}

\define \Ž {\'e}
\define \ {\`e}
\define \ˆ {\`a}
\define \ {\`u}

\define \™ {\^o}
\define \ {\^e}
\define \š {\"o}

\UseAMSsymbols
\topmatter
\title Ramsey, for Auld Lang Syne \endtitle

\abstract\nofrills A stroll taken around the landscape of
Ramsey's Theory. One way, \lq\lq Down from infinite to
finite", then, another way, \lq\lq Up from disorder to
order". An expos\Ž \ made at the
\lq\lq Rencontres arithm\Ž tique et combinatoire", 
Saint-Etienne, june 2006.\endabstract

\author Labib Haddad  \endauthor

\address {120 rue de Charonne, 75011 Paris, France; e-mail:
labib.haddad\@wanadoo.fr } \endaddress

\endtopmatter
\document
\refstyle{A}
\widestnumber\key{ABCD}

\

\heading Ramsey's theorems \endheading

\

\subheading{1 Introduction} 

\

At any given time, in a given
meeting, some of the pairs of persons have already shaken
hands, others have not. Pick any group of six persons. You
are sure to find one of the two following (not
exclusive) situations : Either three of them have
already shaken hands together, or three of them have not.
Stated otherwise, either there is a trio who have or a trio
who have not. Here is the explanation. Concentrate on one
member
$M$ of the group.  Then either $M$ has shaken hands with at
least three  of the five others, or he has not with three
of them. Say he has shaken hands with $P$, $Q$ and $R$. If
one of the three pairs, say $P$ and $Q$, have shaken
hands together, we are done, with the trio $M,P,Q,$ who have.
Otherwise
$P,Q,R,$ is our trio who have not!

\

This is a well-known popular teaser, and an easy way to
introduce what is known as \lq\lq Ramsey's Theorem".

\

In fact, Ramsey proved \bf two \rm theorems of this kind
(see [9]), one in a finite setting , the other in the
infinite setting . More about that in a moment.

\

So six persons is enough to get trios. A natural question
arises : How about quartets? How many to be sure to get
quartets?

\

Let us elaborate. For shortness sake, define a \bf
$2$-coloured graph
\rm to be a \sl complete \rm symmetric graph whose edges are
coloured either red or blue. Of course, a \bf monochromatic
graph
\rm is one all of whose edges have the same colour. Also
define the
\bf size
\rm of a graph to be the number of its vertices. So, in any
$2$-coloured graph of size six, one is sure to find a
monochromatic triangle. That is the essence of our teaser.
The simplest version of Ramsey's Theorem goes like this.

\

\subheading{2 Theorem} Given any integer $h > 0$, there is
an integer $k > 0$ such that each $2$-coloured graph of size
at least $k$ contains a monochromatic subgraph of size
$h$. 

\

The least of all such integers $k$ is denoted $N(2,h)$ and
called the Ramsey number (for $2$-coloured graphs, relative
to
$h$). So, $N(2,3) \leq 6$ and, in
fact, it is easily seen that
$N(2,3)$ is equal to $6$. It is also known that
$N(2,4) = 18$, for instance. But no \sl simple \rm way is
known which would determine the value of $N(2,h)$ as a
function of $h$. The hunt for Ramsey numbers looks somewhat
like that for prime numbers. Compared to the quantity of
primes which we know, that of Ramsey numbers already \lq\lq
caught" is tiny, and almost no useful criteria are available.

\

A \sl generalisation-oriented \rm mind will certainly ask why
should one stick to two colours? In fact, the theorem
extends to $r$-coloured graphs as well. An edge in a
(ordinary) graph is a sort of \sl link \rm that binds a pair
of vertices. Again, why stick to those links by pairs. One
can as well think of another sort of link, by bundles of $n$
vertices. The theorem is still valid in this new setting as
we will readily see.

\

\subheading{3 A bit of terminomogy} For brevity's sake, and
simplicity, we shall deviate again from the traditional
vocabulary and notations. An
\bf $r$-coloured $n$-graph \rm is a configuration defined by
a set
$V$ (the \sl vertices\rm), together with the set $\CP_n(V)$
of all the $n$-element subsets of $V$ (the \sl edges\rm), and
a covering of $\CP_n(V)$ by $r$ sets
$C_1,\dots,C_r$ (the \sl colours\rm) which might
possibly overlap (that is, an edge might have more than one
colour), and even overflow. The size of an
$n$-graph is, again, the number of its vertices, be it
finite or not. Of course, a
\bf monochromatic \rm $n$-graph
is one all of whose edges have a same colour (that is, it has
$\CP_n(V) \subset C_i$ for some index $i$).

\

So, a $2$-coloured graph is nothing else but a $2$-coloured
$2$-graph (provided some of the edges might be, at the same
time, red and blue.)

\

One word more about \bf subgraphs\rm. Take
any subset $W \subset V$ of the set of vertices $V$ of a
given
$r$-coloured $n$-graph $G$. \bf Restricting \rm the set of
vertices to
$W$ means considering the $r$-coloured $n$-graph whose set
of vertices is $W$ and whose colours are
still $C_1,C_2,\dots,C_r$. The new
$r$-colored $n$-graph thus obtained will be called the
subgraph on $V$ \bf restricted \rm to $W$ or, equally well,
the subgraph \bf spanned \rm by $W$ in $G$. Those are all of
what will be called
\bf subgraphs
\rm in this context.

\

The full Ramsey Theorem now reads like this.

\

\subheading{4 Ramsey's Theorem. The Finite Version} Given
integers $r,n,h,$ there is an integer $k$ such that any
$r$-coloured $n$-graph of size at least
$k$ contains a monochromatic $n$-subgraph of size $h$.

\

The least of all such integers $k$ is denoted $N(r,n,h)$ and
called the Ramsey number (for $r$-coloured $n$-graphs,
relative to
$h$). 

\

[Let us mention, for the record, that a Ramsey Number of the
kind
$f(n,q_1,\dots,q_r)$ has also been introduced to be the
least integer $k$ for which every $r$-coloured $n$-graph of
size $k$ contains a monochromatic $n$-subgraph of size
$q_i$ and colour $C_i$ for at least one of the indices $i$.
Then, of course, the inequality $f(n,q_1,\dots,q_r) \leq
N(r,n,\max
\{q_1,\dots,q_r \})$ holds.]

\

The theorem applies equally well, of course, to $r$-coloured
$1$-graphs. But what is a $1$-graph? Well it is much the
same as a set
$V$ of vertices. Then an $r$-coloured $1$-graph is a set
covered by $r$ sets. What the theorem says about this
configuration is precisely Dedekind's \bf pigeon-hole
principle \rm : Given $r$ drawers and an integer $h$, if $k$
articles are to be distributed in the drawers and if $k >
r(h-1)$, then one of the drawers must contain at least $h$ of
those articles, [\sl le principe des tiroirs \rm being the
French name for this principle]. 

\

\smc H. J. RYSER \rm [10, page 38] puts it this way
: \lq\lq Ramsey's theorem may be regarded as a profound
generalization of this simple principle."

\

Incidentally, we have $N(r,1,h) = r(h-1) + 1$.

\

As for the Infinite Version of Ramsey's Theorem, it
concisely says that
$$N(r,n,\infty) = \infty,$$
which, in expanded form, means the following.

\

\subheading{5 Ramsey' Theorem. The Infinite Version} Each
infinite $r$-coloured $n$-graph contains an infinite
monochromatic $n$-subgraph.

\

\heading Down from infinite to finite \endheading

\

Usually, texbooks on Combinatorial Theory contain proofs of
the Finite Version of Ramsey's Theorem. See, for instance,
\smc L. COMTET
\rm [2], \smc M\rm arshall 
\smc HALL\rm, Jr. [7], and \smc H. J. RYSER \rm [10].
Compare with the short proof for the Infinite Version given
below. 

\

As is often the case, the arithmetic of the infinite
(cardinals) seems simpler to handle than that of the finite
integers. This is no exception. Not only is the statement of
the infinite version of the theorem simpler, but its
proof is also much more straightforward (see the \smc
Appendix \rm below). Moreover, as it turns out, with a bit of
ingenuity, from the Infinite Version, one  obtains an easy
proof of the theorem in the Finite version, \sl with a
little bonus\rm. 

\

\subheading{6 Ramsey's Theorem, strenghened} Given integers
$r,n,h,$ there is an integer $k$ with the following
property :

\

$P(r,n,h;k)$ : Any $r$-coloured
$n$-graph whose vertices are $1,2,\dots,k,$
contains a mono\-chromatic $n$-subgraph spanned by $H \subset
\{1,2,\dots,k\}$ and such that
$$|H| \geq h \ \ \text{and} \ \ |H| \geq \min(H).$$

\

Of course, it does not really matter where the 
vertices come from, as long as they are labelled with the
integers from $1$ to $k$.

\

Here is (more than merely a sketch of) the proof.

\

\subheading{7 The proof} Consider a suitable nonstandard
extension
$^*\N$ of the natural integers $\N$. Any ultrapower of $\N$
relative to an incomplete ultrafilter would do. Let
$\Omega = \ ^*\N \setminus
\N$. Then take any
\bf hyperfinite \rm integer $m \in \Omega$
and an
$r$-coloured $n$-graph whose set of vertices is $V =
\{0,1,\dots,m\}$. The $r$-coloured $n$-subgraph
obtained by restricting the set of vertices to $\N$ is
infinite and thus, by the Infinite Version of Ramsey's
Theorem, contains an infinite monochromatic
$n$-subgraph with vertices in $\N$ :
$$v_1 < v_2 < \dots < v_i < \dots.$$
Take $q$ to be a (finite) integer larger than both $h$ and
$v_1$. Then $H = \{v_1,v_2,\dots,v_q\}$ spans a
monochromatic $n$-subgraph with
$$|H| \geq h \ \ \text{and} \ \ |H| \geq \min(H).$$
This simply means that property $P(r,n,h;m)$ is satisfied
by each hyperfinite integer $m \in \Omega$. The following
subset, which is standard, 
$$M = \{ m \in \ ^*\N : P(r,n,h;m) \ \ \text{is
satisfied}\}$$
thus contains all the hyperfinite integers in
$\Omega$. So now, $k = \min (M)$ is the required (standard)
integer : Indeed, $P(r,n,h;k)$ is satisfied and $k$ belongs
to
$\N$ for, otherwise,
$k-1$ would still belong to $\Omega$!\qed

\

\subheading{8 One or two hints more} Let us add some few
words of explanation.

\

Though some familiarity with nonstandard methods
could help in grasping the essence of this (short) proof,
very little
is needed to understand it. Indeed, all the knowledge
needed can be summed up, loosely, as follows.

\

\bf A statement concerning ordinary $\N$ is true for $\N$ if
and only if it is true for nonstandard $^*\N$\rm.

\

\noindent To even better understand the last part of the
proof, it might be helpful to consider it as pertaining to
the
\sl adage \rm : A statement
$S(m)$ true for each infinitely large integer $m$ must be
true for at least one finite integer $k$. In fact, the
saying goes even further : Such a statement as $S(m)$ is true
for all infinitely large integers $m$ if and only if it is
true for all sufficiently large finite integers.

\

\subheading{9 Proofs for the Finite Version} The proofs
usually given for the Finite Version of Ramsey's Theorem do
not use its Infinite Version but rely mainly on recursion,
going from $r$-coloured to
$(r+1)$-coloured and from
$n$-graphs to $(n+1)$-graphs, along the
combinatorial way. Those proofs are said to be
\sl elementary\rm, with different acceptances for the term
\lq\lq elementary". For instance, one such proof can be given
using nothing more than \bf PA\rm, Peano Arithmetic (not Pure
Arithmetic!), which is the first order theory of arithmetic,
built on the axioms of Peano.

\

But, as strange as it may seem, the Strenghened Version of
the Theorem is unprovable in
\bf PA\rm. This undecidability result was established by \smc
J. Paris
\rm and
\smc L. Harrington \rm [8]. 

\

Let us stress this fact again: Ramsey's is a theorem
in \bf PA\rm, while the statement of its strenghtened
version is undecidable in
\bf PA\rm. This also means that Ramsey's Theorem is true
in every extension of $\N$, while its strenghened version 
is true only in some of the extensions and false in others
[not of the kind of ultrapowers, of course].

\

Also notice the following : True, G\š del proved the
incompleteness of \bf PA
\rm in the thirtees (and much more, of course), but the
statements he used for that purpose did not look exactly like
\sl ordinary \rm mathematical statements. To many, they
still may seem a little bit
\sl artificial \rm or too much \bf ad hoc\rm, somehow. Paris
and Harrington, for the first time (to my knowledge) have
shown that a quite ordinary mathematical statement can be
undecidable in \bf PA\rm. It is one of the very first
incompleteness results for
\bf PA \rm which produces a \sl natural \rm undecidable
statement. 

\

\subheading{10 Whither undecidability}

\

The layman will, inevitably, wonder how such statements can
be proved unprovable! Well, to figure it out, just a
glimpse, think of the arithmetical functions that can be
handled by
\bf PA\rm. They are infinite in number, no doubt, but still
denumerable. So, they do not cover the whole \sl spectrum
\rm of possible arithmetical functions. Using Cantor's
Diagonal Argument, one can define, outside \bf PA\rm, a
function that grows faster than all of them. Thus, from the
outside, \sl so to say\rm, one can push integers up, or down,
further than any function from the inside would do!

\

G\š del himself uses Cantor's Diagonal Argument for his
proofs which, nevertheless, are of a different kind, more
syntactical than functional.

\

As Arithmeticians have long known, it is often
easier to prove a theorem of arithmetic using alien
tools, such as analytic tools for instance. Now, everybody
should be aware that some of those theorems do indeed \bf
require
\rm such powerful tools which cannot be dispensed with.

\

Although Peano Arithmetic can serve the purpose of building
a large part of mathematics as it goes, it does not even
cover all of Arithmetic. There is a hierarchy in the affairs
of proofs which, taken at its height, leads on to fortune.
Zermelo-Frankel's set theory
\bf ZFC
\rm is a much
\sl stronger \rm theory than
\bf PA\rm, even though it is not the top, of course. But
\lq\lq this is another story".

\

\heading Up from disorder to order \endheading

\

We have just seen an instance of the \bf downgoing \rm
influence of the Infinite on the Finite. Let me tell you now
about an example of an \bf upgoing \rm construction leading
from disorder to order.

\

\subheading{11 Ramsey and the ordinals} Among the very many
generalizations of Ramsey's Theory, one is about ordinals.
The $n$-graphs in this generalization have their vertices
\bf well-ordered\rm, instead of just being labelled by
integers. To make a long story short, given an ordinal
$\gamma$, let us define an
$(n,\gamma)$-graph to be an $n$-graph with a well-order
of type $\gamma$ on the vertices. Of course, this can be
identified with an $n$-graph whose vertices \bf are \rm the
ordinals $\kappa \in
[0,\gamma[$ from $0$ to $\gamma$ not included.

\

\smc P. Erd\š s \rm and \smc R. Rado \rm
[3,
Corollary to Theorem 39], offered the following
generalization.

\

\subheading{12 Ramsey for ordinals} Given ordinals
$\alpha,\beta$, there is an ordinal $\gamma$ such that each
$2$-coloured $(n,\gamma)$-graph either contains
a red
$(n,\alpha)$-subgraph or a blue
$(n,\beta)$-subgraph.

\

Sticking to the case where $n = 2$, define  
$$f(\alpha,\beta) \ \ \text{to be the least of all
ordinals} \ \  \gamma \ \ \text{such that each
$2$-coloured $(2,\gamma)$-graph}$$ 
\ \ \ \ \ \ \ \ \ \ \ \ \ \ \ \ either contains a red
$(2,\alpha)$-subgraph or a blue $(2,\beta)$-subgraph. 

\

[This
is an analog of Ramsey Numbers $f(2,q_1,q_2)$ mentionned
above.]

\

\subheading{13 An example} Since thirty seven years now, in a
joint note,
\smc L. Haddad
\rm and
\smc G. Sabbagh \rm [5,(1969)], there is a half-page proof of
the fact that
$$f(m,\omega^2) =
\omega^2 \ \ \text{for each integer} \ \ m > 0.$$
As it was, the proof was written
down, bare, with no comments or hints. Our note got very
little attention, in fact it got almost none! So no further
details were ever published. True, the result was already
known : \smc E. Specker \rm [11] had
already given a proof of it, using  an ultrafilter, among
other ingredients. 

\

Here is an \sl expanded \rm form of this short proof in [5],
hoping it will thus illustrate a construction, going
upwards, from disorder to order.

\

\subheading{14 $\bold {\slanted f}(\bold
m, \boldsymbol\omega^2) \boldkey =
\boldsymbol \omega^2$} 

\

Consider a
$2$-coloured $(2,\omega^2)$-graph whose set of vertices is
$\N \times \N$, ordered lexicographically, and an integer $m
> 0$. We are bound to show the following : Either there is a
red subgraph of size $m$ or else a
blue $(2,\omega^2)$-sugraph. We shall show a bit more.

\

Either there are red subgraphs of each size $m > 0$ or else
a blue $(2,\omega^2)$-subgraph.

\

We cover the set $\CP_4(\N)$ of all quadruples with \bf
eight \rm subsets in the following way : Let
$$X = \{n_1,n_2,n_3,n_4\} \ , \ \ n_1 <
n_2 < n_3 < n_4.$$
Starting
from $(n_1,n_2,n_3,n_4)$, one gets $4! = 24$ permutations.
We keep an eye on three of them :
$$X_1 = (n_1,n_2,n_3,n_4) \ , \ X_2 = (n_1,n_3,n_2,n_4) \ ,
\ X_3 = (n_1,n_4,n_2,n_3).$$ 
Let then $X_i = (a,b,a',b')$ for some $i \in \{1,2,3\}$. Two
things can occur : Either the edge $((a,b),(a',b'))$ is red,
and we set $\overline {X_i} = +$ [for red], or it is blue,
and we set
$\overline {X_i} = -$ [for blue]. We call the ordered triple
$(\overline{X_1},\overline{X_2},\overline{X_3})$ the \bf
signature \rm of $X$ [a kind of multicolour]. Since there are
only
$2^3 = 8$ possible signatures, we get a covering of
$\CP_4(\N)$ with those $8$ \sl multicolours \rm. Let
us ponder a moment on the
$8$-coloured $4$-graph we thus obtain and whose set of
vertices is $\N$.

\

What does it really mean that $X = \{n_1 < n_2 < n_3 < n_4\}$
has signature $(+,-,+)$ [that is, (red,blue,red)],
for instance. Well, this means that the edge
$$((n_1,n_2),(n_3,n_4)) \  \text{is red} \ ,
\ ((n_1,n_3),(n_2,n_4) \  \text{is blue} \ ,
\ ((n_1,n_4),(n_2,n_3)) \  \text{is red}.$$
We have \sl transfered\rm, in a sense, (part of) the
structure of the initial $2$-coloured $(2,\omega^2)$-graph on
$\N \times
\N$ to an $8$-coloured $4$-graph on $\N$. We now use Ramsey's
Theorem to get an infinite subset $A$ of $\N$ all of whose
quadruples have a same signature, say $(c_1,c_2,c_3)$.
Restricting the initial graph to the set of vertices $A
\times A$, a moment's thought is enough to see that we can
suppose we already started with a
$2$-coloured $(2,\omega^2)$-graph such that all of
$\CP_4(\N)$ have that same signature $(c_1,c_2,c_3)$.
Now, only two cases can occur : 

\

\noindent 1) Either there is an index $i$ such that $c_i =
+$. 

\

\noindent 2) Or, for each $i$, $c_i = -$.

\

So, the proof has two more steps. In the first case, we
show that there are red subgraphs of each (finite) size
$m > 0$. In the second case, we show that, unless there
are red subgraphs of each (finite) size $m > 0$, there must
be a blue
$(2,\omega^2)$-subgraph.

\

\noindent \bf The first case\rm. Suppose $c_1 = +$. That is,
all the edges $((a,b),(a',b'))$ such that $a < b < a' < b',$
are red. So the \bf infinite \rm sugraph with vertices
$(0,1),(2,3),(3,4),\dots,(2k,2k+1),\dots$ is red.

\

\noindent Suppose $c_2 = +$. That is, all the edges
$((a,b),(a',b'))$ such that $a  < a' < b < b',$ are red.
The subgraph with vertices
$(0,m),(1,m+1),\dots,(k,m+k),\dots,(m-1,2m-1),$ is red and of
size $m$.

\

\noindent Suppose $c_3 = +$. This means that the edges
$((a,b),(a',b'))$ such that $a  < a' < b' < b,$ are red.
This time, the subgraph with vertices
$(0,2m),(1,2m-1),\dots,(k,2m-k),\dots,(m-1,m+1),$ is red of
size $m$.

\

\noindent \bf The second case\rm. Suppose $c_1 = c_2 = c_3 =
-$. This means that, whenever  
$$a < b < a' < b' \ \ \text{or} \ \ a  < a' < b < b' \ \
\text{or} \ \  a  < a' < b' < b,$$
then the edge $((a,b),(a',b'))$ is blue. In order to settle
matters easily, we restrict the vertices of the graph to the
subset
$$L = \{(p,p^n) : \ p \ \ \text{prime and} \ \ n > 1\}.$$
Suppose that the size of red subgraphs does not exceed $m$.
Consider the subgraphs $G(p)$ restricted to $L(p) = \{
(p,p^n)
\ :
\ n > 1\}$. Each one of them contains an infinite blue
subgraph, according to Ramsey's Theorem. So, there is an
infinite subset $M(p)$ of $L(p)$  such that the restriction
of $G(p)$ to $M(p)$ is an infinite blue subgraph.
Set $M = \bigcup M(p)$. This is a subset of $L$
and is well-ordered of type $\omega^2$. 

\

\noindent All we have to do now
is to show that the restriction of the graph on $L$ to $M$ 
is blue. Take any two different vertices $(p,p^s) \neq
(q,q^t)$ in $M$. If
$p = q$, then both vertices are in $M(p)$ and the edge $e =
((p,p^s),(q,q^t))$ is blue. Otherwise, let $p < q$. Then
either 
$$p < p^s < q < q^t \ \ \text{or} \ \ p < q < p^s < q^t \ \
\text{or} \ \ p < q < q^t < p^s$$
holds : In all three cases
the edge $e$ is blue, due to the signature.\qed

\

\heading Appendix \endheading

\

Here is a proof of the Infinite Version of
Ramsey' Theorem :

\

 $R(r,n)$ Each
infinite $r$-coloured $n$-graph contains an infinite
monochromatic $n$-subgraph.

\

\subheading{15 Proof} [The result is obvious for $r = 1$ and
each $n > 0$ because a $1$-coloured $n$-graph is already
monochromatic!]

\

For $n= 1$ and each $r > 0$, the vertices are
covered by $r$ subsets one of which must be infinite, so the
result obtains. We then proceed by induction on $n$.

\

Suppose the result is proved for given $r$ and $n$. Take any
infinite $r$-coloured $(n+1)$-graph whose set of vertices is
$V$ and colours are $C_1,C_2,\dots,C_r$. Starting from
any vertex $v_0$, and using
$R(r,n)$, define by induction
a sequence of vertices 
$$v_0,v_1,\dots,v_k,\dots,$$ 
a sequence of infinite subsets 
$$V \supset V_0 \supset V_1 \supset \dots \supset V_k \supset
\dots,$$
and a function $c : \N \to \{1,\dots,r\}$ such that 
$v_k \in V_{k-1} \setminus V_k$ and
the
$(n+1)$-element subset
$\{v_k\}
\cup A$ has colour
$C_{c(k)}$  for each $A \in \CP_n(V_k)$. At least one of the
subsets $M(i) = \{ k \in \N : c(k) = i \}$ is infinite, say
$M(j)$. Set $W = \{v_k : k \in M(j)\}$ :  This is an
infinite subset of $V$ and the sugraph restricted to $W$ is
clearly monochromatic with colour
$C_j$.\qed

\

\

This proof is exercice \noo 28 in \smc Bourbaki \rm [1, E
III.92]. It should be compared to the proof suggested by \smc
Bourbaki \rm [1, E III.86, exercice 17] for the Finite
Version.

\

\

\centerline{\bf Here is a French version of the text\rm}

\

\

\centerline{\bf Ramsey, for Auld Lang Syne \rm (French
version)}

\

\head R\'esum\'e\endhead 

\

Une promenade dans le d\Ž cor des
th\Ž or\ mes de Ramsey, en descendant, de \lq\lq l'infini au
fini", puis en remontant, du \lq\lq d\Ž sordre vers
l'ordre". On y d\Ž veloppe, entre autre, une tr\ s ancienne
d\Ž monstration du fait que $f(m,\omega^2) =
\omega^2$ pour tout entier $m > 0$.

\

\heading Pour un expos\Ž \ aux \lq\lq Rencontres arithm\Ž
tique et combinatoire" de
\ Saint-Etienne, juin 2006\endheading

\heading une version en fran\c cais \endheading

\

\heading Les Th\Ž or\ mes de Ramsey \endheading

\

\subheading{1 Introduction}

\

A un moment donn\Ž , dans une assembl\Ž e quelconque,
quelques unes des paires de ses membres ont d\Ž j\ˆ \ fait
connaissance et d'autres pas. Choisissons au hasard six de
ces membres. On est alors assur\Ž \ de se trouver
dans l'une ou l'autre des deux situations
suivantes (voire les deux \ˆ \ la fois) : ou bien
trois de ces six membres ont d\Ž j\ˆ \ fait connaissance
entre eux ou bien trois d'entre eux ne l'ont pas encore
fait. Autrement dit, ou bien il y a un trio qui a fait
connaissance ou bien un trio qui ne l'a pas fait. Voici
l'explication. Portons notre attention sur l'un de ces six
membres, soit
$M$. Ou bien $M$ a d\Ž j\ˆ \
fait la connaissance d'au moins trois des cinq autres ou
bien il ne l'a pas fait avec trois autres. Supposons qu'il
ait d\Ž j\ˆ \ fait connaissance avec $P$, $Q$ et $R$. Si deux
d'entre ces derniers, disons
$P$ et $Q$, ont d\Ž j\ˆ \
fait connaissance, nous tenons un trio $M,P,Q,$ qui l'ont d\Ž
j\ˆ
\ fait. Sinon, $P,Q,R,$ est un trio qui ne l'ont pas fait!

\

C'est un de ces \sl petits probl\ mes \rm assez
r\Ž pandu et bien connu. C'est aussi un moyen rapide pour
faire conna\" tre ce que l'on appelle, commun\Ž ment,
\lq\lq le th\Ž or\ me de Ramsey".

\

En r\Ž alit\Ž , Ramsey a \Ž tablit \bf deux \rm th\Ž or\
mes de ce genre (voir [9]), l'un dans le cadre du fini,
l'autre dans le cadre infini. On en reparlera un plus tard.

\

Ainsi six membres suffisent pour avoir des trios. Une
question vient naturellement \ˆ \ l'esprit : qu'en est-il
des quatuors? Combien de membres pour \ tre certain d'avoir
des quatuors?

\

Rentrons dans les d\Ž tails. Pour faire court, on appellera
\bf graphe bicolore \rm tout graphe \sl complet \rm et sym\Ž
trique dont les ar\ tes sont color\Ž es en rouge ou en
bleu. Bien entendu, un graphe \bf monochrome \rm est un
graphe dont toutes les ar\ tes sont d'une m\ me couleur.
La \bf taille \rm d'un graphe d\Ž signe le nombre de ses
sommets. Ainsi, dans un graphe bicolore quelconque de taille
six, on est s\ž r de trouver un triangle monochrome. C'est
le contenu de ce petit probl\ me. La version la plus simple
du th\Ž or\ me de Ramsey dit ceci.

\

\subheading{2 Th\Ž or\ me} Pour tout entier $h > 0$, il
existe un entier $k > 0$ tel que tout graphe bicolore de
taille (au moins) \Ž gale \ˆ \
$k$ contient un sous-graphe monochrome de taille $h$.

\

Le plus petit de tous ces entiers $k$ est d\Ž sign\Ž \ par
$N(2,h)$ et s'appelle le nombre de Ramsey (pour les graphes
bicolores, relatif \ˆ \ $h$). Ainsi, $N(2,3) \leq 6$ et on
voit facilement que
$N(2,3)$ est en fait \Ž gal \ˆ \ $6$. On sait \Ž galement
que $N(2,4) = 18$, par exemple. Cependant, on ne conna\" t
encore aucun moyen \sl simple \rm pour d\Ž terminer la
valeur de $N(2,h)$ en fonction de $h$. La chasse aux nombres
de Ramsey ressemble un peu \ˆ \ celle des nombres premiers.
Cela \Ž tant, compar\Ž e \ˆ \ la quantit\Ž \ des nombres
premiers connus, celle des nombres de Ramsey d\Ž j\ˆ \
\lq\lq captur\Ž s" est
\sl infime\rm, et l'on ne dispose quasiment d'aucun crit\
re utile. 

\

Un esprit tourn\Ž \ vers la g\Ž n\Ž ralisation (plut\™ t que
les sp\Ž cialisations) ne peut s'emp\ cher de se demander
pourquoi on devrait s'en tenir
\ˆ \ deux couleurs seulement. Le th\Ž or\ me s'\Ž tend tout
aussi bien aux graphes $r$-colores. Une ar\ te dans un
graphe (ordinaire) est une sorte de \sl lien \rm entre une
paire de sommets. De m\ me, pourquoi s'en tenir
\ˆ \ ces liens par paires. On peut tr\ s bien envisager une
autre sorte de lien, par paquets de
$n$ sommets. Le th\Ž or\ me vaut \Ž galement dans ce cadre
g\Ž n\Ž ral comme on va le voir tout de suite.

\

\subheading{3 Un brin de terminologie} Pour demeurer bref, et
rester simple, on va s'\Ž carter de nouveau des notations et
du vocabulaire traditionnels. Un
\bf $n$-graphe $r$-colore \rm est une configuration d\Ž
finie par la donn\Ž e d'un ensemble $S$ (les \sl
sommets\rm), de l'ensemble $\CP_n(S)$ des parties de $S$ de
cardinal
$n$ (les \sl ar\ tes\rm), ainsi que d'un recouvrement de
$\CP_n(S)$ par $r$ ensembles $C_1,\dots,C_r$ (les \sl
couleurs\rm) qui peuvent chevaucher (de
sorte qu'une m\ me ar\ te puisse avoir, \Ž ventuellement, 
plusieurs couleurs
\ˆ \ la fois) et m\ me d\Ž border. La taille d'un
$n$-graphe d\Ž signe toujours le nombre de ses sommets,
qu'il soit fini ou infini. Bien entendu, un \bf $n$-graphe
monocolore \rm est celui dont toutes les ar\ tes ont une
m\ me couleur (autrement dit, tel que $\CP_n(S) \subset
C_i$ pour un indice $i$ donn\Ž ).

\

Ainsi, un graphe bicolore n'est rien autre qu'un $2$-graphe
$2$-colore (en acceptant qu'une ar\ te puisse \ tre, tout
\ˆ
\ la fois, rouge et bleue).

\

Un mot encore au sujet des \bf sous-graphes\rm. Soit $T
\subset S$ une partie quelconque de l'ensemble $S$ des
sommets d'un $n$-graphe $r$-colore $G$. \bf Restreindre \rm
l'ensemble des sommets \ˆ \ $T$ veut dire considerer le
$n$-graphe $r$-colore dont l'ensemble des sommets est $T$ et
dont les couleurs sont toujours $C_1,C_2,\dots,C_r$. Le
nouveau $n$-graphe
$r$-colore ainsi obtenu sera applel\Ž \ le sous-graphe sur
$S$ \bf restreint \rm \ˆ \
$T$, ou encore, le sous-graphe \bf sous-tendu \rm par
$T$ dans $G$. Ce sont les seuls qui seront appel\Ž s
\bf sous-graphes \rm dans ce contexte.

\

Le th\Ž or\ me de Ramsey en sa g\Ž n\Ž ralit\Ž \ se pr\Ž
sente comme suit.

\

\subheading{4 Le th\Ž or\ me de Ramsey. La version \lq\lq
finie"}
Les entiers
$n,r,h,$ \Ž tant donn\Ž s, il existe un entier $k$ tel que
tout $n$-graphe
$r$-colore de taille $k$ contienne un $n$-sous-graphe
monochrome de taille $h$.

\

Le plus petit de ces entiers $k$ est d\Ž sign\Ž \ par
$N(n,r,h)$ et se nomme le nombre de Ramsey (pour les
$n$-graphes $r$-colores, relatif \ˆ \ $h$).

\

[On mentionnera, pour m\Ž moire, qu'un nombre de Ramsey
du genre
$f(n,q_1,\dots,q_r)$ a \Ž t\Ž \ \Ž galement introduit pour
d\Ž signer le plus petit entier $k$ tel que tout 
$n$-graphe $r$-colore de taille $k$ poss\ de un
$n$-sous-graphe monochrome de taille $q_i$ et de couleur
$C_i$ pour l'un au moins des indices $i$. Bien entendu, on a
$f(n,q_1,\dots,q_r) \leq N(n,r,\max \{q_1,\dots,q_r \})$.]

\

\'Evidemment, le th\Ž or\ me s'applique tout aussi bien aux
$1$-graphes $r$-colores. Mais qu'est-ce
qu'un $1$-graphe? Eh bien! c'est, \ˆ \ tout prendre, la m\
me chose qu'un ensemble $S$ de sommets. Ainsi, un $1$-graphe
$r$-colore est un ensemble recouvert par $r$ ensembles.
Ce que le th\Ž or\ me dit \ˆ \ propos de cette configuration
est pr\Ž cis\Ž ment \bf le principe des tiroirs \rm de
Dedekind : \Ž tant donn\Ž \ un colombier ayant $r$ nids, et
un entier
$h$, si $k$ pigeons habitent ce colombier
et si
$k > r(h-1)$, alors l'un de ces nids devra abriter au
moins
$h$ pigeons. [En anglais, ce principe porte le nom de \sl
pigeon-hole principle\rm.]

\

Signalons ce mot de \smc H. J. RYSER \rm [10, page 38] : 
le th\Ž or\ me de Ramsey peut \ tre consid\Ž r\Ž \ comme
une g\Ž n\Ž ralisation profonde de ce principe simple. 

\

Incidemment, on a $N(1,r,h) = r(h-1) + 1$.

\

Pour ce qui est du th\Ž or\ me de Ramsey en sa version
\lq\lq infinie", il dit succintement que l'on a
$$N(n,r,\infty) = \infty,$$
ce qui, sous une forme d\Ž velopp\Ž e, veut dire ceci.

\

\subheading{5 Le th\Ž or\ me de Ramsey. La version \lq\lq
infinie"} Tout $n$-graphe $r$-colore infini contient un
$n$-sous-graphe monochrome infini.

\

\heading Descente de l'infini au fini \endheading

\

Habituellement, les manuels d'analyse combinatoire
contiennent des d\Ž monstrations de la version "finie" du
th\Ž or\ me de Ramsey. Voir, par exemple, \smc L. COMTET
\rm [2], \smc M\rm arshall 
\smc HALL\rm, Jr. [7] et \smc H. J. RYSER \rm [10]. On pourra
comparer ces d\Ž monstrations avec celle, courte, de la
version "infinie" pr\Ž sent\Ž e ci-dessous.

\

Comme c'est souvent le cas, l'arithm\Ž tique de l'infini
(pour les cardinaux) para\" t plus simple \ˆ \ manier que
celle des entiers finis. Le cas pr\Ž sent ne fait pas
exception. Non seulement l'\Ž nonc\Ž \ de la version \lq\lq
infinie" du th\Ž or\ me est plus simple, mais sa d\Ž
monstration est \Ž galement plus directe (voir l'\smc
Appendice \rm ci-dessous). De plus, avec un peu d'ing\Ž
niosit\Ž , il se trouve que, de la version
\lq\lq infinie" du th\Ž or\ me, on tire une d\Ž
monstration facile du th\Ž or\ me dans sa version \lq\lq
finie", avec \sl une petite prime \rm en plus.

\

\subheading{6 Le th\Ž or\ me de Ramsey
renforc\Ž } Des entiers
$n,r,h,$ \Ž tant donn\Ž s, il existe un entier
$k$ ayant la propri\Ž t\Ž \ suivante :

\

$P(n,r,h;k)$ : Tout
$n$-graphe $r$-colore dont les sommets sont 
$1,2,\dots,k,$ contient un $n$-sous-graphe monochrome
sous-tendu par $H \subset
\{1,2,\dots,k\}$ et tel que
$$|H| \geq h \ \ \text{et} \ \ |H| \geq \min(H).$$

\

Bien entendu la provenance des sommets importe peu, il
suffit qu'ils soient num\Ž rot\Ž s \ˆ \ l'aide des entiers
de $1$ \ˆ \ $k$.

\

Voici (davantage qu'une simple esquisse de) cette d\Ž
monstration.

\

\subheading{7 La d\Ž monstration} On consid\ re une
extension convenable $^*\N$ des entiers naturels $\N$. Une
ultrapuissance quelconque de $\N$ relativement \ˆ \ un
ultrafiltre incomplet peut faire l'affaire. On pose $\Omega
= \ ^*\N \setminus \N$. On prend ensuite un entier
\bf hyperfini \rm quelconque
$m \in \Omega$ et un $n$-graphe $r$-colore dont l'ensemble
des sommets est $S = \{0,1,\dots,m\}$. Le $n$-sous-graphe
$r$-colore obtenu par restriction de l'ensemble des sommets
\ˆ \
$\N$ est infini de sorte que, d'apr\ s la version \lq\lq
infinie" du th\Ž or\ me de Ramsey, il contient un
$n$-sous-graphe monochrome infini ayant (dans $\N$) les
sommets 
$$s_1 < s_2 < \dots < s_i < \dots.$$
On prend un entier (fini) $q$ plus grand que $h$ et $s_1$.
Ainsi $H = \{s_1,s_2,\dots,s_q\}$ sous-tend un
$n$-sous-graphe monochrome et l'on a
$$|H| \geq h \ \ \text{et} \ \ |H| \geq \min(H).$$
Cela veut simplement dire que la propri\Ž t\Ž \ $P(n,r,h;m)$
est satisfaite par chacun des entiers hyperfini $m \in
\Omega$. Le sous-ensemble (standard)
$$M = \{ m \in \ ^*\N : P(r,n,h;m) \ \ \text{est
satisfaite}\}$$ contient ainsi tous les entiers hyperfinis
de $\Omega$.  De sorte que $k = \min (M)$ est l'entier
(standard) annonc\Ž \ : en effet, $P(n,r,h;k)$ est
satisfaite et $k$ appartient bien \ˆ \ $\N$ car, sinon, $k-1$
appartiendrait toujours \ˆ \ $\Omega$!\qed

\

\subheading{8 Une ou deux indications
compl\Ž mentaires} Ajoutons quelques mots d'explication.

\

Bien qu'une certaine familiarit\Ž \ avec les m\Ž thodes
nonstandards puisse aider \ˆ \ mieux saisir l'essence de
cette (courte) d\Ž monstration, tr\ s peu est requis pour la
comprendre. En effet, tout ce que l'on a besoin de savoir
peut se r\Ž sumer, un peu rapidement, comme suit.

\

\bf Un \Ž nonc\Ž \ relatif \ˆ \ $\N$ est vrai pour $\N$
si et seulement s'il est vrai pour $^*\N$
nonstandard.\rm

\

\noindent Afin de mieux comprendre encore la derni\ re
partie de la d\Ž monstration, il peut \ tre utile de la
consid\Ž rer comme relevant de l'\sl adage \rm : un \Ž
nonc\Ž \ $E(m)$ vrai pour tout entier $m$ infiniment grand
doit \ tre vrai pour au moins un entier $k$ fini. En r\Ž
alit\Ž , le proverbe va plus loin : un \Ž nonc\Ž \ tel que
$E(m)$ est vrai pour tous les entiers infiniment grands  si
et seulement s'il est vrai pour les entiers finis
suffisamment grands.

\

\subheading{9 Les d\Ž monstrations de la version \lq\lq
finie"} Les d\Ž monstrations que l'on donne habituellement
de la version \lq\lq finie" du th\Ž or\ me de Ramsey
n'utilisent pas sa version \lq\lq infinie", mais reposent
surtout sur la r\Ž currence, allant des $n$-graphes aux
$(n+1)$-graphes et de $r$ couleurs en $(r+1)$ couleurs,
\sl combinatoirement\rm. Ces preuves sont dites \sl \Ž l\Ž
mentaires \rm en diverses acceptions du mot \lq\lq \Ž l\Ž
mentaire". En particulier, une telle preuve existe qui
n'utilise rien de plus que l'arithm\Ž tique de Peano, dite
\bf PA\rm, qui est la th\Ž orie du premier ordre de
l'arithm\Ž tique, construite \ˆ \ l'aide des axiomes de
Peano.

\

Cependant, aussi \Ž trange que cela puisse para\" tre, la
version renforc\Ž e du th\Ž or\ me n'est pas d\Ž montrable
dans \bf PA\rm. Cette ind\Ž cidabilit\Ž \ a \Ž t\Ž \ \Ž
tablie par \smc J. Paris \rm et \smc L. Harrington \rm [8].

\

Insistons encore sur ce fait : le th\Ž or\ me
de Ramsey est un th\Ž or\ me de \bf PA \rm tandis que
l'\Ž nonc\Ž \ de sa version renforc\Ž e est ind\Ž cidable
dans
\bf PA\rm. Cela veut \Ž galement dire que le th\Ž or\ me de
Ramsey est vrai dans toute extension de $\N$, tandis que sa
version renforc\Ž e est vraie dans certaines extensions et
fausse dans d'autres [pas du genre
ultrapuissances, bien entendu].

\

On observera \Ž galement ceci : il est vrai que G\š del
a \Ž tabli l'incompl\Ž tude de \bf
PA \rm  (et bien davantage encore,
sans aucun doute,) dans les ann\Ž es trente, mais les \Ž
nonc\Ž s qu'il a utilis\Ž s pour ce faire ne ressemblent pas
tout
\ˆ \ fait \ˆ \ des \Ž nonc\Ž s math\Ž matiques \sl
ordinaires\rm. Pour beaucoup, ils peuvent encore appara\"
tre un peu \sl artificiels \rm ou quelque peu trop \bf ad
hoc\rm. Paris et Harrignton, pour la premi\ re fois
(pour autant que je le sache) ont montr\Ž
\ qu'un \Ž nonc\Ž \ math\Ž matique tout \ˆ \ fait ordinaire
peut \ tre ind\Ž cidable dans \bf PA\rm.

\

\subheading{10 D'o\  \ vient l'ind\Ž cidabilit\Ž }

\ 

Le profane s'\Ž tonnera, in\Ž vitablement : comment peut-on
d\Ž montrer que de tels \Ž nonc\Ž s sont ind\Ž
montrables. Eh bien! pour en avoir une id\Ž e, un tr\ s
rapide coup d'\oe il, on doit penser aux fonctions arithm\Ž
tiques que
\bf PA
\rm peut manier. Elles sont en nombre infini, certainement,
mais toujours d\Ž nombrables. Elles ne couvrent donc pas
tout le
\sl spectre \rm des fonctions arithm\Ž tiques possibles. En
utilisant le proc\Ž d\Ž \ diagonal de Cantor, on peut d\Ž
finir, \ˆ \ l'ext\Ž rieur de \bf PA\rm, une fonction qui
cro\" t plus vite que chacune d'elles. Ainsi, de l'ext\Ž
rieur, \sl pour ainsi dire\rm, on peut pousser des entiers
vers le haut, ou vers le bas, plus loin qu'aucune des
fonctions de l'int\Ž rieur ne le ferait!

\

G\š del lui-m\ me utilise une forme de proc\Ž d\Ž \
diagonal pour ses d\Ž monstrations lesquelles sont, n\Ž
anmoins, d'une autre nature, davantage syntaxique que
fonctionnelle.

\

Les arithm\Ž ticiens le savent depuis longtemps, il est
souvent plus facile de d\Ž montrer un th\Ž or\ me en
utilisant des outils \Ž trangers, tels que des outils
analytiques, par exemple. Maintenant tout le monde devrait
se rendre compte que quelques uns de ces th\Ž or\ mes n\Ž
cessitent de tels outils puissants dont on ne peut pas se
dispenser.

\

Bien que l'arithm\Ž tique de Peano puisse servir le projet
de b\‰ tir une grande partie des math\Ž matiques telles
qu'elles sont, elle ne couvre m\ me pas toute l'arithm\Ž
tique. En mati\ re de preuves, il y a une hierarchie qui,
dans les sommets, conduit au succ\ s. La th\Ž orie des
ensembles \bf ZFC \rm de Zermelo-Frankel est bien plus forte
que \bf PA\rm, m\ me si ce n'est pas le sommet, bien s\ž r.
Mais \lq\lq ceci est une autre histoire".

\

\heading Remont\Ž e du d\Ž sordre vers l'ordre \endheading

\

Nous venons de voir un exemple de l'influence \bf
descendante \rm de l'infini sur le fini. Voici, \ˆ \ pr\Ž
sent, l'exemple d'une construction qui \bf remonte \rm du
d\Ž sordre vers l'ordre.

\

\subheading{11 Ramsey et les ordinaux} Parmi les tr\ s
nombreuses g\Ž n\Ž ralisations de la th\Ž orie de Ramsey, il
en est une relative aux ordinaux. Dans cette g\Ž ralisation,
les $n$-graphes ont leurs sommets \bf bien ordonn\Ž s \rm,
au lieu d'\ tre seulement num\Ž rot\Ž s \ˆ \ l'aide
d'entiers. Pour faire court, \Ž tant donn\Ž \ un ordinal
$\gamma$, convenons d'appeler $(n,\gamma)$-graphes les
$n$-graphes dont les sommets sont bien ordonn\Ž s suivant
l'ordre
$\gamma$. Bien entendu, ceux-ci peuvent \ tre identifi\Ž s
aux $n$-graphes dont les sommets \bf sont \rm les ordinaux
$\kappa \in [0,\gamma[$ de $0$ \ˆ \ $\gamma$ non inclus.

\

\smc P. Erd\š s \rm et \smc R. Rado \rm
[3, corollaire au th\Ž or\ me 39], ont pr\Ž sent\Ž \ la g\Ž
n\Ž ralisation suivante.

\

\subheading{12 Ramsey et les ordinaux} \'Etant donn\Ž s des
ordinaux $\alpha,\beta$, il existe un ordinal $\gamma$ tel
que tout $(n,\gamma)$-graphe 
bicolore contienne ou bien un
$(n,\alpha)$-sous-graphe rouge ou bien un
$(n,\beta)$-sous-graphe bleu.

\

En s'en tenant au cas o\ \ $n = 2$, on d\Ž finit
$$f(\alpha,\beta) \ \ \text{le plus petit des
ordinaux} \ \  \gamma \ \ \text{tels que chaque
$(2,\gamma)$-graphe bicolore}$$ 
\ \ \ \ \ \    contient ou bien un
$(2,\alpha)$-sous-graphe rouge ou bien un
$(2,\beta)$-sous-graphe bleu. 

\

[C'est l'analogue des nombres
de Ramsey $f(2,q_1,q_2)$ mentionn\Ž s ci-dessus.]

\

\subheading{13 Un exemple} Depuis trente sept ans, dans une
note conjointe,
\smc L. Haddad
\rm et
\smc G. Sabbagh \rm [5,(1969)], on trouve une d\Ž monstration
d'une demi-page qui \Ž tablit le fait que l'on a
$$f(m,\omega^2) =
\omega^2 \ \ \text{pour tout entier} \ \ m > 0.$$
Telle quelle, la d\Ž monstration apparaissait nue, r\Ž
dig\Ž e sans commentaires ni indications. Notre note n'attira
que tr\ s peu d'attention, en r\Ž alit\Ž , elle n'en obtint
presque aucune. Aussi, aucun autre d\Ž tail ne fut-il jamais
publi\Ž . Il est vrai que ce r\Ž sultat \Ž tait d\Ž j\ˆ \
connu: \smc E. Specker \rm [11] en
avait donn\Ž \ une d\Ž monstration, faisant intervenir un
ultrafiltre, parmi d'autres ingr\Ž dients.

\

Voici une forme \sl d\Ž ploy\Ž e \rm de cette courte d\Ž
monstration dans [5], en esp\Ž rant ainsi qu'elle
puisse illustrer une construction remontant du d\Ž sordre
vers l'ordre.

\

\subheading{14 $\bold {\slanted f}(\bold
m, \boldsymbol\omega^2) \boldkey =
\boldsymbol \omega^2$} 

\

On consid\ re un $(2,\omega^2)$-graphe bicolore dont
l'ensemble des sommets est $\N \times \N$, ordonn\Ž \
lexicographiquement, et un entier $m > 0$. Nous devons
montrer ceci : il existe un sous-graphe rouge de
taille $m$, sinon il existe un $(2,\omega^2)$-sous-graphe
bleu. Nous allons \Ž tablir un petit peu plus.

\

Il y a des sous-graphes rouges de toute taille finie
$m > 0$, sinon il existe un $(2,\omega^2)$-sous-graphe bleu.

\

On recouvre l'ensemble $\CP_4(\N)$ de tous les  quadruplets
\ˆ \ l'aide de \bf huit \rm sous-ensembles de la mani\ re
suivante. Soit
$$X = \{n_1,n_2,n_3,n_4\} \ , \  n_1 <
n_2 < n_3 < n_4.$$
En commen\c cant par $(n_1,n_2,n_3,n_4)$, on
obtient $4! = 24$ permutations. On garde un \oe il sur trois
d'entre elles :
$$X_1 = (n_1,n_2,n_3,n_4) \ , \ X_2 = (n_1,n_3,n_2,n_4) \ ,
\ X_3 = (n_1,n_4,n_2,n_3).$$
Soit alors $X_i = (a,b,a',b')$ pour un indice $i \in
\{1,2,3\}$ donn\Ž . Deux choses peuvent arriver: ou bien
l'ar\ te $((a,b),(a',b'))$ est rouge, et on posera
$\overline {X_i} = +$ [pour le rouge], ou bien elle est bleu,
et on posera $\overline {X_i} = -$ [pour le bleu]. On appelle
le triplet $(\overline{X_1},\overline{X_2},\overline{X_3})$
la
\bf signature \rm de $X$ [une esp\ ce de multicouleur].
Comme il n'y a que $2^3 = 8$ signatures possibles, on
obtient un recouvrement de $\CP_4(\N)$ \ˆ \ l'aide de ces $8$
multicouleurs. Scrutons un peu le $4$-graphe
$8$-colore que l'on obtient et dont l'ensemble des sommets
est $\N$.

\

Que signifie exactement le fait que $X = \{n_1 < n_2 < n_3 <
n_4\}$ ait pour signature $(+,-,+)$ [autrement
dit, (rouge,bleu,rouge)], par exemple. Et bien! cela veut
dire que l'ar\ te
$$((n_1,n_2),(n_3,n_4)) \  \text{est rouge} \ ,
\ ((n_1,n_3),(n_2,n_4) \  \text{est bleu} \ ,
\ ((n_1,n_4),(n_2,n_3)) \  \text{est rouge}.$$
Nous avons \bf transf\Ž r\Ž  \rm, en un sens, (une partie
de) la structure du graphe de d\Ž part sur $\N \times \N$ \ˆ
\ un
$4$-graphe
$8$-colore sur $\N$. Utilisons, \ˆ \ pr\Ž sent, le th\Ž or\
me de Ramsey afin d'obtenir une partie infinie $A$ de $\N$
dont tous les quadruplets ont une m\ me signature, soit
$(c_1,c_2,c_3)$. En restreigant le graphe initial
\ˆ \ l'ensemble des sommets $A \times A$, un moment de r\Ž
flexion suffit pour voir que l'on peut supposer
\ tre parti, d'embl\Ž e, d'un $(2,\omega^2)$-graphe
bicolore pour lequel l'ensemble tout entier $\CP_4(\N)$ a la
m\ me signature
$(c_1,c_2,c_3)$. \`A pr\Ž sent, seul deux cas peuvent se pr\Ž
senter :

\

\noindent 1) Ou bien il y a un indice $i$ pour lequel $c_i =
+$. 

\

\noindent 2) Ou bien on a $c_i = -$ pour chaque $i$.

\

Il faut donc encore deux \Ž tapes dans la d\Ž monstration.
Dans le premier cas, on montre qu'il y a des
sous-graphes rouges de toute taille (finie) $m > 0$. Dans
le second cas, on montre qu'il existe n\Ž cessairement un
$(2,\omega^2)$-sous-graphe bleu, \ˆ \ moins qu'il n'y ait
des sous-graphes rouges de toute taille (finie) $m > 0$.

\

\noindent \bf Le premier cas\rm. On suppose que $c_1 = +$. 
Cela veut dire que les ar\ tes $((a,b),(a',b'))$
pour lesquelles $a < b < a' < b',$ sont rouges. Ainsi, le
sous-graphe \bf infini \rm dont les sommets sont
$(0,1),(2,3),(3,4),\dots,(2k,2k+1),\dots$ est-il rouge.

\

\noindent On suppose que $c_2 = +$. Cela veut dire que
toutes les ar\ tes $((a,b),(a',b'))$ pour lesquelles  $a  <
a' < b < b',$ sont rouges. Le sous-graphe dont les sommets
sont $(0,m),(1,m+1),\dots,(k,m+k),\dots,(m-1,2m-1),$ est donc
rouge et de taille $m$.

\

\noindent On suppose que $c_3 = +$. Cela veut dire que les
ar\ tes $((a,b),(a',b'))$ pour lesquelles $a  < a' < b' <
b,$ sont rouges. Cette fois, le sous-graphe dont les sommets
sont
$(0,2m),(1,2m-1),\dots,(k,2m-k),\dots,(m-1,m+1),$ est rouge
et de taille $m$.

\

\noindent \bf Le second cas\rm. On suppose que $c_1 = c_2 =
c_3 = -$. Cela veut dire que, d\ s que l'on a
$$a < b < a' < b' \ \ \text{ou} \ \ a  < a' < b < b' \ \
\text{ou} \ \  a  < a' < b' < b,$$
l'ar\ te  $((a,b),(a',b'))$ est bleue. Pour r\Ž gler les
choses simplement, on restreint les sommets du graphe au
sous-ensemble 
$$L = \{(p,p^n) : \ p \ \ \text{premier et} \  n > 1\}.$$
On suppose que la taille des sous-graphes rouges ne d\Ž
passe pas un entier donn\Ž \ $m$. On consid\ re le
sous-graphe
$G(p)$ restreint \ˆ \ $L(p) = \{ (p,p^n)
\ :
\ n > 1\}$. Chacun de ces graphes contient un sous-graphe
infini bleu, d'apr\ s le th\Ž or\ me de Ramsey. Il existe
donc une partie infinie $M(p)$ de $L(p)$ telle que la
restriction de $G(p)$ \ˆ \ $M(p)$ soit un sous-graphe infini
bleu. Soit
$M = \bigcup M(p)$. C'est une partie de $L$ et elle est bien
ordonn\Ž e de type $\omega^2$.

\

\noindent Tout ce qui nous reste \ˆ \ faire c'est de d\Ž
montrer que la restriction du graphe sur $L$ \ˆ \ $M$ est
bleue. Soient $(p,p^s) \neq
(q,q^t)$ deux sommets distincts dans $M$. Si $p = q$, les
deux sommets sont dans $M(p)$ et l'ar\ te
$u =((p,p^s),(q,q^t))$ est donc
bleue.  Sinon, soit
$p < q$. Alors, ou bien on a
$$p < p^s < q < q^t \ \ \text{ou} \ \ p < q < p^s < q^t \ \
\text{ou} \ \ p < q < q^t < p^s,$$
et, dans les trois cas, l'ar\ te $u$ est bleue, gr\‰ ce \ˆ \
la signature.\qed

\

\heading Appendice \endheading

\

Voici une d\Ž monstration de la version \lq\lq infinie" du
th\Ž or\ me de Ramsey :

\

$R(n,r)$ Tout $n$-graphe $r$-colore infini contient un
$n$-sous-graphe monochrome infini.

\

\subheading{15 D\Ž monstration} [Le r\Ž sultat est \Ž vident
pour $r = 1$ quel que soit $n > 0$ car un $n$-graphe
monocolore est monochrome!]

\

Pour $n = 1$ et chaque $r > 0$, les sommets sont recouverts
par $r$ sous-ensembles dont l'un au moins doit \ tre
infini, ce qui donne le r\Ž sultat. On proc\ de ensuite par
r\Ž currence sur $n$.

\

On suppose le r\Ž sultat acquis pour $n$ et $r$ donn\Ž s.  On
se donne un $(n+1)$-graphe  $r$-colore quelconque dont
l'ensemble des sommets est $S$ et les couleurs sont
$C_1,C_2,\dots,C_r$. On commence par un sommet quelconque
$s_0$ et, en utilisant
$R(n,r)$, on d\Ž finit par r\Ž currence une suite de sommets
$$s_0,s_1,\dots,s_k,\dots$$ 
une suite de parties infinies
$$S \supset S_0 \supset S_1 \supset \dots \supset S_k
\supset \dots$$
et une application $c : \N \to
\{1,\dots,r\}$ tels que $s_k \in S_{k-1} \setminus S_k$ et,
pour chaque $A \in \CP_n(V_k)$, la partie
$\{s_k\} \cup A$ de cadinal
$(n+1)$ ait la couleur $C_{c(k)}$. Au moins l'une des
parties $M(i) = \{ k \in \N : c(k) = i \}$ est infinie, soit
$M(j)$. On pose $T = \{v_k : k \in M(j)\}$ : c'est une
partie infinie de $S$ et le sous-graphe restreint \ˆ \ $T$
est visiblement monochrome, de couleur $C_j$.\qed

\

Cette d\Ž monstration est l'exercice \noo 28 dans \smc
Bourbaki \rm [1, E III.92]. On devrait la comparer avec la
d\Ž monstration sugg\Ž r\Ž e par \smc
Bourbaki \rm [1, E III.86, exercice 17] pour la version
\lq\lq finie".

\

\enddocument